\documentclass[11pt]{amsart}
\usepackage{amsmath,amssymb,latexsym,soul}
\usepackage[initials]{amsrefs}
\usepackage{fancyvrb}
\usepackage[left=2.6cm,right=2.6cm,top=2.7cm,bottom=2.7cm]{geometry}

\usepackage{color,enumitem,graphicx}
\usepackage[colorlinks=true,urlcolor=blue,
citecolor=red,linkcolor=blue,linktocpage,pdfpagelabels,bookmarksnumbered,bookmarksopen]{hyperref}
\usepackage[english]{babel}

\usepackage[hyperpageref]{backref}

\numberwithin{equation}{section}
\newtheorem{theorem}{Theorem}[section]

\newcommand{\eps}{\varepsilon}

\newcommand{\Ra}{R_{\alpha,p,q}}
\newcommand{\Sao}{S_{\alpha,p,q,1}}
\newcommand{\Sam}{S_{\alpha,p,q,m}}
\newcommand{\San}{S_{\alpha,p,q,n}}
\newcommand{\Sar}{S_{\alpha,p,q,\infty}}
\newcommand{\Szo}{S_{0,p,q,1}}
\newcommand{\I}{{\mathbb I}}

\newcommand{\R}{{\mathbb R}}
\newcommand{\N}{{\mathbb N}}
\newcommand{\Int}{\mathrm{Int}}

\title[Symmetry of $n$-mode positive solutions for  
2D H\'enon systems]{Symmetry of $n$-mode positive solutions \\
for two-dimensional H\'enon type systems}

\author[N.\ Shioji]{Naoki Shioji}
\address{Department of Mathematics
\newline\indent
Faculty of Engineering
\newline\indent
Yokohama National University
\newline\indent
Tokiwadai, Hodogaya-ku, Yokohama 240-8501, Japan}
\email{shioji@math.sci.ynu.ac.jp}

\author[M.\ Squassina]{Marco Squassina}
\address{Dipartimento di Informatica
\newline\indent
Universit\`a degli Studi di Verona
\newline\indent
C\'a Vignal 2, Strada Le Grazie 15, 37134 Verona, Italy}
\email{marco.squassina@univr.it}

\thanks{The first author was partially supported by 
the Grant-in-Aid for Scientific Research (C) (No.\
21540214) from Japan Society for the Promotion of Science.
The second author was partially supported by 2009 italian MIUR project:
   ``Variational and Topological Methods in the Study of Nonlinear Phenomena''}
\subjclass[2000]{35B06, 35B09, 35J15}
\keywords{Radial symmetry, H\'enon type systems, $n$-mode positive solutions}

\begin{document}

\begin{abstract}
We provide a symmetry result for $n$-mode positive solutions of a general class of semi-linear
elliptic systems under cooperative conditions on the nonlinearities. Moreover, we apply the result
to a class of H\'enon systems and provide the existence of multiple $n$-mode positive solutions.
\end{abstract}

\maketitle

\section{Introduction and results}

\noindent
Let $D=\{x \in \R^2: |x|<1\}$ and, for $m\in \N$ with $m\geq 2$, consider the system
\begin{equation}
\label{problem}
\left\{
\begin{aligned}
\Delta u^i + f^i(|x|,u^i) & =0 && \text{in $D$,}\\
u^i&=0 && \text{on $\partial D$,}
\end{aligned}
\right. \,\,\,\quad\text{for $i=1,\ldots,m$,}
\end{equation}
where $f^i$ are smooth functions over $(0,1)\times (0,\infty)^m$.\ Semi-linear elliptic systems as \eqref{problem} arise
naturally in many physical and biological contests, see 
e.g.~\cites{murray,smoller,Rad1,Rad2,Rad3} and the references therein.
As far as the symmetry of positive solutions is concerned and the functions $f^i$ are decreasing in the radial variable, the
celebrated moving plane method \cite{GNN} can be applied when the system
is cooperative namely $\partial f^i/\partial u^j\geq 0$ 
for every  $i\neq j$ \cites{Def1,troy,KesPac}. 
The aim of this note is to establish a general symmetry result (Theorem~\ref{tone})
for $n$-mode ($2\pi/n$-rotation invariant) solutions, namely solutions $(u^1,\dots,u^m)$
such that each component $u^i:\overline{D}\to\R$, in polar coordinates, satisfies
$$
u^i(r,\theta)=u^i(r,\theta+2\pi/n), \,\,\quad \text{for all $(r,\theta)\in [0,1]\times\R$,}
$$
as well as provide a meaningful application of it (Theorem~\ref{applicationHenon})
to the system of H\'enon type 
\begin{equation}\label{eq:1.1}
\left\{
\begin{aligned}
\Delta u+\frac{2p}{p+q}|x|^\alpha u^{p-1}v^q&=0 && \text{in $D$,}\\
\Delta v+\frac{2q}{p+q}|x|^\alpha u^pv^{q-1}&=0 && \text{in $D$,}\\
u>0,\,\, v&>0 && \text{in $D$,} \\
u=v&=0 && \text{on $\partial D$.}
\end{aligned}
\right.
\end{equation}
Quite recently, these systems were carefully investigated 
in~\cites{yanfu,WY} (see also \cites{alves,han1,Han2} and 
references therein) and they can be considered as a vectorial
counterpart of the celebrated equation 
\begin{equation*}
\left\{
\begin{aligned}
\Delta u +|x|^\alpha u^{p-1}&=0 && \text{in $D$,}\\
u&>0 && \text{in $D$,} \\
u&=0  && \text{on $\partial D$,}
\end{aligned}
\right.
\end{equation*}
first studied in \cite{Ni} after being introduced by
H\'enon in \cite{henon} in connection with the research of rotating stellar structures.
We shall say that $u$ is of class $C^n$ at the origin if $u$ is of
class $C^{n-1}$ in a neighborhood of the origin and each
$(n-1)$-th partial derivative is
totally differentiable at the origin. Then we prove the following

\noindent
\begin{theorem}
\label{tone}
Let $m,n \in \N$ with $m,n\geq 2$ and
$f^1,\ldots,f^m \in C((0,1)\times (0,\infty)^m,\R)$ 
such that
\begin{enumerate}
\item 
for each $i\in \{1,\ldots,m\}$ and $(u^1,\ldots,u^m) \in (0,\infty)^m$, the map
$$
r\mapsto r^{2-2n}f^i(r,u^1,\ldots,u^m): (0,1)\rightarrow \R
$$
is nonincreasing;
\item
for each $i\in \{1,\ldots,m\}$ and $r \in (0,1)$,
$f^i(r,\cdot,\ldots,\cdot)\in C^1((0,\infty)^m,\R)$;
\item
for each 
$i,j\in \{1,\ldots,m\}$ with $i\neq j$
and $(r,u^1,\ldots,u^m) \in (0,1)\times (0,\infty)^m$,
$$
\frac{\partial f^i}{\partial u^j}(r,u^1,\ldots,u^m)\geq 0;
$$
\item
for each $i,j\in\{1,\ldots,m\}$, $r_0 \in (0,1)$ and $M\in (0,\infty)$,
$$
\sup\Big\{\Big|\frac{\partial f^i}{\partial u^j}
(r,u^1,\ldots,u^m)\Big|:
(r,u^1,\ldots,u^m) \in (r_0,1)\times (0,M]^m \Big\}<\infty.
$$
\end{enumerate}
Let $(u^1,\ldots,u^m) \in C^2(D\setminus\{0\})\cap C(\overline{D})$
be a solution of \eqref{problem}
such that each $u^i$ is $n$-mode, positive and 
of class $C^n$ at the origin. Then, each
$u^i$ is radially symmetric and $\frac{\partial u^i}{\partial r}(|x|)<0$ for $r=|x|$.
\end{theorem}

\noindent
For scalar equations, this result was obtained in \cite{SW}. Due to the recent interest
of the community for the symmetry issues for elliptic systems, we believe that the statement above 
is of interest. Also, it admits
some interesting consequences, see for instance Theorem~\ref{applicationHenon} below.
Of course, system~\eqref{problem} includes both variational and nonvariational problems or systems
of Hamiltonian type, see e.g.\ \cite{Def-Mon} for a wide overview. We point out, in particular, that the weakly coupled 
semi-linear Schr\"odinger systems, see \cite{mps} and the references therein, 
which come from physically relevant situations and have recently received much attention, satisfy 
conditions (ii)-(iv). 

For the sake of completeness, we refer the reader to \cites{DGP1,DGP2,DP} for recent partial (foliated Schwarz symmetry) symmetry results
for the smooth solutions to \eqref{problem} in rotationally invariant domains and for possibly sign-changing 
solutions and where the maps $r\mapsto f^i(r,s_1,\dots,s^m)$ are possibly nondecreasing and some convexity
assumptions are assumed on the $s_i$ variables.
\vskip5pt
\noindent
For every $\alpha\geq 0$ and $p,q>1$, let us set
$$
\Ra(u,v):=\frac{\displaystyle\int_D (|\nabla u|^2+|\nabla v|^2)}
{\Big(\displaystyle\int_D|x|^\alpha|u|^p|v|^q\Big)^\frac{2}{p+q}},
\,\,\,\quad\text{for any $u,v\in H^1_0(D)$.}
$$
Moreover, for each $\gamma>0$, we shall denote by $\lceil \gamma\rceil$
the smallest integer greater than or equal to $\gamma$.
Then, we have the following

\begin{theorem}
\label{applicationHenon}
The following facts hold.
\begin{enumerate}
 \item[$(\I)$] 
If $\alpha\in (0,\infty)$, $p,q\in (1,\infty)$
and $(u,v)$ is an $n$-mode solution of \eqref{eq:1.1}
with $n\geq 1+\lceil \alpha/2\rceil$,
then $u,v$ are radially symmetric.
\item[$(\I\I)$]  
For each $\alpha\in (2,\infty)$ and $p,q\in (1,\infty)$,
if $n_\alpha \geq 1$
then \eqref{eq:1.1} has a nonradial
$n$-mode 
solution $(u_n,v_n)$ for $n=1,\ldots,n_\alpha$ such that
\begin{equation}
\label{order}
\Ra(u_1,v_1)< \cdots<\Ra(u_{n_\alpha},v_{n_\alpha}),
\end{equation}
where 
$n_\alpha$ is the greatest integer less than
\begin{equation}
\label{m-alpha}
\Big(\frac{\alpha+2}{2\alpha}\Big)^\frac{4}{p+q-2}
\Big(\frac{\alpha-2}{\alpha}\Big)^\frac{2\alpha}{p+q-2}
\Big(1+\frac{\alpha}{2}\Big).
\end{equation}
In particular, the following facts hold.
\begin{enumerate}
\renewcommand{\labelenumii}{(\roman{enumii})}
\renewcommand{\theenumii}{\roman{enumii}}
\item 
For each $\alpha\in (2,\infty)$,
if at least one of $p,q\in (1,\infty)$ is large enough,
then $n_\alpha=\lceil \alpha/2\rceil$,
that is \eqref{eq:1.1} has a nonradial 
$n$-mode solution
for $n=1,\ldots,
\lceil \alpha/2\rceil$
satisfying \eqref{order}.
\item
For each $p,q \in (1,\infty)$ and $\ell \in \N$,
if $\alpha \in (2,\infty)$ is large enough
\eqref{eq:1.1} has a nonradial
$n$-mode solution $(u_n,v_n)$
for $n=1,\ldots,\ell$ such that 
$$
\Ra(u_1,v_1)< \cdots<\Ra(u_\ell,v_\ell).
$$
In particular,
for each $p,q\in (1,\infty)$,
the number of nonradial solutions of \eqref{eq:1.1}
tends to infinity as $\alpha\rightarrow \infty$.
\end{enumerate}
\end{enumerate}
\end{theorem}

\noindent
Hence, as far as $\alpha$ gets large, the symmetry breaking phenomenon occurs
and we can find as many $n$-mode positive solutions as we want.
In Section~\ref{provathm1} we shall prove Theorem~\ref{tone},
while in Section~\ref{provathm2} we shall provide the proof of 
Theorem~\ref{applicationHenon}.

\section{Proof of Theorem~\ref{tone}}
\label{provathm1}
\noindent
By using the planar polar coordinates,
for each $i=1,\ldots,m$,
we define the function $\widetilde{u}^i:\overline{D}\rightarrow \R$
by setting $\widetilde{u}^i(r,\theta):=u^i(r^{1/n},\theta/n)$,
for every $(r,\theta)\in \overline{D}$.
Since $(u^1,\ldots,u^m)$ satisfies system~\eqref{problem}
and each $u^i$ is $n$-mode,
we can see that
$\widetilde{u}^i \in C^2(D\setminus\{0\})\cap C(\overline{D})$
and $(\widetilde{u}^1,\ldots,\widetilde{u}^m)$ satisfies the system
\begin{equation}
\label{new-problem}
\left\{
\begin{aligned}
\Delta \widetilde{u}^i + \widetilde{f}^i(|x|,\widetilde{u}^1,\ldots,
\widetilde{u}^m) & =0 && 
\text{in $D\setminus\{0\}$,}\\
\widetilde{u}^i&=0 && \text{on $\partial D$,}
\end{aligned}
\right.\qquad\text{for $i=1,\ldots,m$.}
\end{equation}
Here,
$\widetilde{f}^i\in C((0,1)\times (0,\infty)^m,\R)$ is the function defined by
$$
\widetilde{f}^i(r, t^1,\ldots,t^m)
:=n^{-2}r^{(2-2n)/n}f^i(r^{1/n},t^1,\ldots,t^m),
$$ 
for every
$(r,t^1,\ldots,t^m)\in
(0,1)\times (0,\infty)^m$
and $\widetilde{f}^i$ satisfies
\begin{equation}
\label{new-decreasing}
\text{for each $(t^1,\ldots,t^m)\in (0,\infty)^m$,
$r\mapsto\widetilde{f}^i(r,t^1,\ldots,t^m)$
is nonincreasing.}
\end{equation}      
Indeed, from
\begin{align*}
\Delta &\widetilde{u}^i(r,\theta) + 
\widetilde{f}^i(r,\widetilde{u}^1(r,\theta),\ldots,
\widetilde{u}^m(r,\theta))\\
& = \frac{1}{n^2}r^\frac{2-2n}{n}
\biggl(
u^i_{rr}\Bigl(\,r^\frac{1}{n},\frac{\theta}{n}\,\Bigr)
+
\frac{1}{r^\frac{1}{n}}
u^i_r
\Bigl(\,r^\frac{1}{n},\frac{\theta}{n}\,\Bigr)
+
\frac{1}{r^\frac{2}{n}}
u^i_{\theta\theta}\Bigl(\,r^\frac{1}{n},\frac{\theta}{n}\,\Bigr)\\
&\qquad\qquad\qquad\qquad
+f^i\Bigl(r^\frac{1}{n},u^1\Bigl(\,r^\frac{1}{n},\frac{\theta}{n}\,\Bigr),
\ldots, u^m\Bigl(\,r^\frac{1}{n},\frac{\theta}{n}\,\Bigr)
\Bigr)
\biggr)
=0,
\end{align*}
we deduce \eqref{new-problem},
and we can easily see that \eqref{new-decreasing} holds as well, in light of assumption \thetag{i}. 
For each $\lambda \in (0,1)$,
we set $\Sigma_\lambda=\{x \in D: x_1>\lambda\}$
and we define the map $h_\lambda:\overline{\Sigma_\lambda}\rightarrow
\overline{D}$
by
$h_\lambda(x)=(2\lambda-x_1,x_2)$ 
for $x=(x_1,x_2)\in \overline{\Sigma_\lambda}$.
We note that $h_\lambda$ satisfies
\begin{equation}
\label{E-norm}
|h_\lambda(x)|< |x| \quad 
\text{for each $\lambda \in (0,1)$
and $x \in \Sigma_\lambda
\cup \Int_{\partial D}(\overline{\Sigma_\lambda}\cap\partial D)$.}
\end{equation}
Here, for a subset $E$ of $\partial D$,
we denote by $\Int_{\partial D}E$,
the interior set of $E$ with respect to the relative 
topology of $\partial D$.
We set 
$x_\lambda=\left(2\lambda,0\right)$ 
for $\lambda \in (0,1)$. We can see
\begin{equation}
\label{some-properties-one}
x_\lambda \in 
\begin{cases}
\Sigma_\lambda & \text{for each $\lambda \in
      (0,\frac{1}{2})$,}\\
\partial \Sigma_{\lambda} & \text{for $\lambda=\frac{1}{2}$,}
\\
\R^2\setminus\overline{\Sigma_\lambda}  &
\text{for each $\lambda \in (\frac{1}{2},1)$}
\end{cases}
\end{equation}
and
\begin{equation}
\label{some-properties-two}
h_\lambda(x_\lambda)=0 \quad
\text{for each $\lambda \in (0,1/2]$.}
\end{equation}
For the sake of completeness,
we note that
$\Sigma_\lambda \setminus\{x_\lambda\}=\Sigma_\lambda$
for each $\lambda \in [\frac{1}{2},1)$
and
$\overline{\Sigma_\lambda} \setminus\{x_\lambda\}
=\overline{\Sigma_\lambda}$
for each $\lambda \in (\frac{1}{2},1)$.
For each $i,j=1,\ldots,m$, 
we define
$v^i_\lambda \in C^2(\Sigma_\lambda\setminus\{x_\lambda\})
\cap C(\overline{\Sigma_\lambda})$
and $c^{ij}_\lambda\in L^\infty(\Sigma_\lambda)$ by setting
\begin{equation}
\label{v}
v^i_\lambda(x):= \widetilde{u}^i(x)-\widetilde{u}^i(h_\lambda(x)),
\,\,\quad\text{for $x \in \overline{\Sigma_\lambda}$,}
\end{equation}
and
\begin{equation}
\label{c}
c^{ij}_\lambda(x)
:=
-\int_0^1 \frac{\partial\widetilde{f}^i}{\partial u^j}
\bigl(
|x|,s\widetilde{u}^1(x)+(1-s)\widetilde{u}^1(h_\lambda(x)),\ldots,
s\widetilde{u}^m(x)+(1-s)\widetilde{u}^m(h_\lambda(x))\bigr)\,ds.
\end{equation}
By the assumptions of Theorem~\ref{tone},
we can see that $c_\lambda^{ij}\leq 0$ if $i\neq j$ for $x\in \Sigma_\lambda$ and
\begin{equation}
\label{supremum}
\sup_{r<\lambda<1}\sup_{x \in
\Sigma_\lambda}|c^{ij}_\lambda(x)|<\infty,
\,\,\quad\text{for each $r\in (0,1)$ and $i,j=1,\ldots,m$.}
\end{equation}
Therefore, it holds
\begin{equation}
\label{diff-inequality}
-\Delta v^i_\lambda(x)+\sum_{j=1}^m c^{ij}_\lambda(x)v^j_\lambda(x)\leq 0 \quad 
\text{for $\lambda \in (0,1)$,
$x \in \Sigma_\lambda\setminus\{x_\lambda\}$
and $i=1,\ldots,m$.}
\end{equation}
Indeed, \eqref{diff-inequality} can be obtained as follows:
{\allowdisplaybreaks
\begin{align*}
0 &=\Delta \widetilde{u}^i(h_\lambda(x))+
\widetilde{f}^i(|h_\lambda(x)|,\widetilde{u}^1(h_\lambda(x)),
\ldots,\widetilde{u}^m(h_\lambda(x))) \\
&\qquad\qquad-\Delta\widetilde{u}^i(x) -
\widetilde{f}^i(|x|,\widetilde{u}^1(x),\ldots,\widetilde{u}^m(x)) \\
&\geq-\Delta v^i_\lambda(x)+
\widetilde{f}^i(|x|,\widetilde{u}^1(h_\lambda(x)),\ldots,
\widetilde{u}^m(h_\lambda(x))) 
-
\widetilde{f}^i(|x|,\widetilde{u}^1(x),\ldots,\widetilde{u}^m(x))\\
&=
-\Delta v^i_\lambda(x)+
\sum_{j=1}^m c^{ij}_\lambda(x)v^j_\lambda(x).
\end{align*}
}
We set
\begin{equation}
\label{Aone-muone}
\begin{aligned}
A_1&= \{ \lambda \in [1/2,1): \text{$v^i_\lambda(x)<0$
for each $x \in \Sigma_\lambda$ and $i\in\{1,\ldots,m\}$}\},\\
\mu_1&=\inf_{\lambda \in A_1}\lambda.
\end{aligned}
\end{equation}
We now claim that $A_1\neq\emptyset$.
Let $i\in \{1,\ldots,m\}$ and
$\lambda \in [1/2,1)$ such that 
$\lambda$ is sufficiently close to $1$.
Then we can easily see
$v^i_\lambda(x)\leq 0$ 
for $x \in \partial\Sigma_\lambda$
and $v^i_\lambda(x)< 0$ for
$x \in \Int_{\partial D}(
\partial D\cap\partial\Sigma_\lambda)$
from \eqref{E-norm}.
Since $|\Sigma_\lambda|\ll 1$ 
and \eqref{diff-inequality}
holds,
by \cite{MR2086750}*{Corollary~14.1}, 
we have $v^i_\lambda\leq 0$ on $\overline{\Sigma_\lambda}$.
By 
\begin{equation}
\label{SMP-inequality}
-\Delta v^i_\lambda(x)+c^{ii}_\lambda(x)v^i_\lambda(x)
\leq -\sum_{j\neq i}^m c^{ij}_\lambda(x)v^j_\lambda(x)\leq 0 
\quad\text{in $\Sigma_\lambda\setminus\{x_\lambda\}$}
\end{equation}
and
the strong maximum principle,
we have $v^i_\lambda<0$ in $\Sigma_\lambda$.
Since $i$ is any element of $\{1,\ldots,m\}$,
we have shown $\lambda \in A_1$, which proves the claim.

\noindent
We now claim that $\mu_1=1/2\in A_1$.
Let $i\in \{1,\ldots,m\}$.
We have $v^i_{\mu_1}(x)\leq 0$ for $x \in \Sigma_{\mu_1}$.
Since
\eqref{SMP-inequality} holds
with $\lambda=\mu_1$
and $v^i_{\mu_1}(x)<0$ for $x \in \Int_{\partial D}(\partial\Sigma_{\mu_1}
\cap \partial D)$ from \eqref{E-norm},
by the strong maximum principle,
we have $v^i_{\mu_1}(x)<0$ for $x \in \Sigma_{\mu_1}$.
Since $i$ is an arbitrary element of $\{1,\ldots,m\}$,
we have $\mu_1\in A_1$.
We will show $\mu_1=1/2$.
Suppose not, namely $\mu_1>1/2$.
Again let $i \in \{1,\ldots,m\}$.
Let $G$ be an open set such that
$\overline{G}\subset \Sigma_{\mu_1}$ 
and $|\Sigma_{\mu_1}\setminus \overline{G}|\ll 1$.
We have
$\max_{x \in \overline{G}}v^i_{\mu_1}(x)<0$.
Let $0<\varepsilon\ll 1$.
Then we have
$\max_{x \in
\overline{G}}v^i_{\mu_1-\varepsilon}(x)<0$
and $|\Sigma_{\mu_1-\varepsilon}\setminus
\overline{G}|\ll 1$.
Since \eqref{diff-inequality} holds with $\lambda=\mu_1-\varepsilon$
and
$v^i_{\mu_1-\varepsilon}(x)\leq 0$ 
for $x \in \partial( \Sigma_{\mu_1-\varepsilon}
\setminus\overline{G})$,
we have
$v^i_{\mu_1-\varepsilon}(x)\leq 0$ 
for $x \in \Sigma_{\mu_1-\varepsilon}$
by 
\cite{MR2086750}*{Corollary~14.1}.
From
$v^i_{\mu_1-\varepsilon}(x)< 0$ 
for $x \in (\Int_{\partial D}(\partial\Sigma_{\mu_1-\varepsilon}
\cap \partial D))\cup\partial G$
and the strong maximum principle,
we have
$v^i_{\mu_1-\varepsilon}(x)<0$ 
for $x \in \Sigma_{\mu_1-\varepsilon}$.
Since $i$ is an arbitrary element of 
$\{1,\ldots,m\}$, we have $\mu_1-\varepsilon \in A_1$.
This is a contradiction. Hence $\mu_1=1/2\in A_1$.

\noindent
We now set 
\begin{equation}
\label{Atwo-mutwo}
\begin{aligned}
A_2&:=\{\lambda 
\in (0,1/2):\;
\text{$v^i_\lambda(x) <0$ 
for each $x \in \Sigma_\lambda$ and $i\in \{1,\ldots,m\}$}\},\\
\mu_2&:=\inf_{\lambda \in A_2}\lambda.
\end{aligned}
\end{equation}
We now claim that $A_2\neq\emptyset$.
Let $i\in \{1,\ldots,m\}$.
We note that $x_{1/2}=(1,0)$.
Let $G$ be an open set such that
$\overline{G}\subset \Sigma_{1/2}$ and
$|\Sigma_{1/2}\setminus G|\ll 1$.
From $1/2\in A_1$ and 
$\overline{G}\subset \Sigma_{1/2}$,
we have $\max_{x \in \overline{G}} v^i_{1/2}(x)<0$.
Let $\lambda \in (0,1/2)$ such that
$\lambda$ is sufficiently close to $1/2$.
We note $|\Sigma_\lambda \setminus\overline{G}|\ll 1$
and $x_\lambda$ is close to $(1,0)$.
We choose a sufficiently small open neighborhood $U$ of $x_\lambda$
with $\overline{U}\subset \Sigma_\lambda$,
and we set $H=G\cup U$.
Then we have 
$v^i_\lambda(x)<0$ for $x \in \overline{H}$,
$v^i_\lambda(x)\leq 0$
for $x \in \partial \Sigma_\lambda\cup \partial H$
and $|\Sigma_\lambda \setminus\overline{H}|\ll 1$.
Since \eqref{diff-inequality} holds on $\Sigma_\lambda\setminus\overline{H}$,
by \cite{MR2086750}*{Corollary~14.1},
we have $v^i_\lambda\leq 0$ on $\Sigma_\lambda$.
From \eqref{SMP-inequality}
and the strong maximum principle,
we have $v^i_\lambda< 0$ on $\Sigma_\lambda$.
Since $i$ is an arbitrary element of $\{1,\ldots,m\}$,
we have shown $\lambda \in A_2$.

\noindent
Recalling that $u$ is of class $C^n$ at the origin,
arguing exactly as in \cite{SW}*{Lemma~4} we get
\begin{equation}
\label{pd}
\frac{\partial (\widetilde{u}^i\circ h_{\mu_2})}{\partial
x_1}(x_{\mu_2})=0
\quad\text{for each $i =1,\ldots,m$.}
\end{equation}

\noindent
We now claim that $\mu_2=0$. Suppose not.
Let $i\in \{1,\ldots,m\}$.
Then 
we have $\mu_2\in (0, 1/2)$
by the previous claim 
and we can see
$v^i_{\mu_2} \leq 0$ 
on $\overline{\Sigma_{\mu_2}}$.
We will show $v^i_{\mu_2}<0$ on $\Sigma_{\mu_2}\setminus\{x_{\mu_2}\}$.
We have $v^i_{\mu_2}(x)<0$ for 
$x \in \Int_{\partial D}(\partial \Sigma_{\mu_2}\cap \partial D)$
from \eqref{E-norm}.
By \eqref{SMP-inequality} with $\lambda=\mu_2$
and the strong maximum principle,
we have $v^i_{\mu_2}<0$ on $\Sigma_{\mu_2}\setminus\{x_{\mu_2}\}$. Next,
we will show
$v^i_{\mu_2}(x_{\mu_2})<0$.
Suppose $v^i_{\mu_2}(x_{\mu_2})<0$ does not hold,
i.e.,
$v^i_{\mu_2}(x_{\mu_2})=0$.
Let $\nu_1=(-1,0)$
and $\nu_2=(1,0)$.
From \eqref{pd},
we have
$$
\begin{aligned}
\frac{\partial v^i_{\mu_2}}{\partial \nu_1}(x_{\mu_2})
&=-\frac{\partial \widetilde{u}^i}{\partial x_1}(x_{\mu_2})
+\frac{\partial (\widetilde{u}^i\circ h_{\mu_2})}{\partial x_1}(x_{\mu_2})
=-\frac{\partial \widetilde{u}^i}{\partial x_1}(x_{\mu_2}),\\
\frac{\partial v^i_{\mu_2}}{\partial \nu_2}(x_{\mu_2})
&=\frac{\partial \widetilde{u}^i}{\partial x_1}(x_{\mu_2})
-\frac{\partial (\widetilde{u}^i\circ h_{\mu_2})}{\partial x_1}(x_{\mu_2})
=\frac{\partial \widetilde{u}^i}{\partial x_1}(x_{\mu_2}).
\end{aligned}
$$
By Hopf's lemma,
we obtain
$$-\frac{\partial \widetilde{u}^i}{\partial x_1}(x_{\mu_2})<0
\quad\text{and}\quad
\frac{\partial \widetilde{u}^i}{\partial x_1}(x_{\mu_2})<0,$$
which is a contradiction.
So we have shown
$v^i_{\mu_2}(x_{\mu_2})<0$.
Thus we have
$v^i_{\mu_2}<0$
on
$\Sigma_{\mu_2}$.
Since $i$ is an arbitrary element of $\{1,\ldots,m\}$,
we have $\mu_2 \in A_2$.
Again let $i \in \{1,\ldots,m\}$.
We choose an open set $G$ such that
$\overline{G}\subset \Sigma_{\mu_2}$ and
$|\Sigma_{\mu_2}\setminus\overline{G}|\ll 1$.
We have
$\max_{\overline{G}} 
v^i_{\mu_2}<0$.
Let $0<\varepsilon\ll 1$.
Then we have
$|\Sigma_{\mu_2-\varepsilon}\setminus\overline{G}| \ll 1$
and
$\max_{\overline{G}}
v^i_{\mu_2-\varepsilon}<0$.
Since \eqref{diff-inequality} holds with $\lambda=\mu_2-\varepsilon$,
by \cite{MR2086750}*{Corollary~14.1}, we have
$v^i_{\mu_2-\varepsilon}(x)\leq 0$
for $x\in
\Sigma_{\mu_2-\varepsilon}\setminus
\overline{G}$.
By \eqref{SMP-inequality} with $\lambda=\mu_2-\varepsilon$
and the strong maximum principle,
we have
$v^i_{\mu_2-\varepsilon}(x)<0$
for $x\in
\Sigma_{\mu_2-\varepsilon}\setminus
\overline{G}$.
Hence we have
shown
$v^i_{\mu_2-\varepsilon}(x)<0$
for
$x \in \Sigma_{\mu_2-\varepsilon}$.
Since $i$ is an arbitrary element of $\{1,\ldots,m\}$,
we have $\mu_2-\varepsilon\in A_2$,
which is a contradiction.
Therefore we obtain $\mu_2=0$.

\noindent
We can finally conclude the proof of Theorem~\ref{tone}.
Let $i \in \{1,\ldots,m\}$.
By the conclusions above,
we can infer that
$\widetilde{u}^i$ is radially symmetric
and $\frac{\partial \widetilde{u}^i}{\partial r}(|x|)<0$ for $r=|x|\in (0,1)$.
From the definition of $\widetilde{u}^i$,
we can find $u^i$ is also radially symmetric and 
$\frac{\partial u^i}{\partial r}(|x|)<0$.

\section{Proof of Theorem~\ref{applicationHenon}}
\label{provathm2}
\noindent
Let us first prove assertion ($\I$) of Theorem~\ref{applicationHenon}.
Assume that $(u,v)$ is an $n$-mode solution to
system~\eqref{eq:1.1} such that $n \geq 1+\lceil \alpha/2
\rceil$. Then, we may choose $\hat n,m\in\N$ such that $m/\hat n\in\N,$ 
$(\alpha+2)\hat n\leq 2n<2(\alpha+2)\hat n$ and
$$
m>\max\Big\{\frac{n\hat n}{(\alpha+2)\hat n-n},\frac{2n}{\alpha+2}\Big\}.
$$
Setting $\hat u(r,\theta):=u(r^{m/n},m\theta/n)$ and $\hat v(r,\theta):=v(r^{m/n},m\theta/n)$, it is readily seen that $\hat u$ and $\hat v$
are both $\hat m=m/\hat n$-mode and solve, in $D\setminus\{0\}$,  the system
\begin{equation*}
\left\{
\begin{aligned}
\Delta \hat u+\frac{2pm^2}{(p+q)n^2}|x|^{\frac{m(\alpha+2)-2n}{n}} \hat u^{p-1}\hat v^q&=0 && \text{in $D\setminus\{0\}$,}\\
\Delta \hat v+\frac{2pm^2}{(p+q)n^2}|x|^{\frac{m(\alpha+2)-2n}{n}} \hat u^p \hat v^{q-1}&=0 && \text{in $D\setminus\{0\}$,}\\
\hat u>0,\,\, \hat v&>0 && \text{in $D\setminus\{0\}$,} \\
\noalign{\vskip2pt}
\hat u=\hat v&=0 && \text{on $\partial D$.}
\end{aligned}
\right.
\end{equation*}
We need to show that $(\hat u,\hat v)$ is a solution of the corresponding system on $D$,
namely
\begin{equation}
\label{systD}
\left\{
\begin{aligned}
\Delta \hat u+\frac{2pm^2}{(p+q)n^2}|x|^{\frac{m(\alpha+2)-2n}{n}} \hat u^{p-1}\hat v^q&=0 && \text{in $D$,}\\
\Delta \hat v+\frac{2pm^2}{(p+q)n^2}|x|^{\frac{m(\alpha+2)-2n}{n}} \hat u^p \hat v^{q-1}&=0 && \text{in $D$.}
\end{aligned}
\right.
\end{equation} 
To this aim, let $\varphi\in C^\infty_c(D)$ a function and let $\eps\in (0,1)$. Then,
if $D_\eps$ denotes the ball centered at zero with radius $\eps$, we get
\begin{align*}
0& =\int_{D\setminus D_\eps}\Delta \hat u\varphi+\frac{2pm^2}{(p+q)n^2}\int_{D\setminus D_\eps}|x|^{\frac{m(\alpha+2)-2n}{n}} \hat u^{p-1}\hat v^q\varphi   \\
& =-\int_{\partial D_\eps}\frac{\partial \hat u}{\partial r}\varphi dS
-\int_{D\setminus D_\eps}\nabla \hat u\nabla\varphi+\frac{2pm^2}{(p+q)n^2}\int_{D\setminus D_\eps}|x|^{\frac{m(\alpha+2)-2n}{n}} \hat u^{p-1}\hat v^q\varphi, \\
0& =\int_{D\setminus D_\eps}\Delta \hat v\varphi+\frac{2pm^2}{(p+q)n^2}\int_{D\setminus D_\eps}|x|^{\frac{m(\alpha+2)-2n}{n}} \hat u^{p}\hat v^{q-1}\varphi   \\
& =-\int_{\partial D_\eps}\frac{\partial \hat v}{\partial r}\varphi dS
-\int_{D\setminus D_\eps}\nabla \hat v\nabla\varphi+\frac{2pm^2}{(p+q)n^2}\int_{D\setminus D_\eps}|x|^{\frac{m(\alpha+2)-2n}{n}} \hat u^{p}\hat v^{q-1}\varphi.
\end{align*}
Since also $u\in C^1(\overline{D})$, the functions
$
\frac{\partial u}{\partial r}(r,\theta),\,
\frac{1}{r}\frac{\partial u}{\partial \theta}(r,\theta),\,
\frac{\partial \varphi}{\partial r}(r,\theta),\,
\frac{1}{r}\frac{\partial \varphi}{\partial \theta}(r,\theta)
$
are bounded on $D$.
From $\bar{u}(r,\theta)=u(r^\frac{m}{n},\frac{m}{n}\theta)$, 
there exists some positive constant $C$ such that
\begin{equation*}
 \Big|\frac{\partial \bar{u}}{\partial r}(r,\theta)\Big|
\leq Cr^{\frac{m}{n}-1},\quad 
 \Big|\frac{\partial \bar{u}}{\partial \theta}(r,\theta)\Big|
\leq Cr,
\end{equation*}
for each $(r,\theta)\in D$.
Hence, we have
\begin{align*}
 \Big|\int_{|x|=\varepsilon}\frac{\partial\bar{u}}{\partial
   r}\varphi\,dS\Big|
&\leq C\varepsilon^\frac{m}{n}, \\
\Big|\int_{\{\varepsilon<|x|<1\}} 
\nabla\bar{u}\,\nabla\varphi
-\int_D
\nabla\bar{u}\,\nabla\varphi
\Big|
& \leq
C\eps^{\frac{m}{n}+1}+C\eps^2.
\end{align*}
Hence, letting $\varepsilon\rightarrow +0$, we conclude that $(\hat u,\hat v)$ is a weak (and hence a strong) solution to \eqref{systD}.
Since $(m(\alpha+2)-2n)/n\geq \hat m-1$, it follows that $\hat u,\hat v\in C^{\hat m}(\overline{D})$. 
Furthermore, since $(\alpha+2)\hat n\leq 2n$ the map 
$r\mapsto r^{2-2\hat m+(m(\alpha+2)-2n)/n}$ is nonincreasing. In turn, by applying 
Theorem~\ref{tone}, it follows that $\hat u$ and $\hat v$ are radially symmetric and hence 
$u$ and $v$ are radially symmetric, concluding the first part of the proof.

\noindent
We now come to the proof of assertion ($\I\I$).\
We set $H_{n}=\{u \in H^1_0(D): \text{$u$ is $n$-mode}\}$
for all $n \in \N$
and $H_{\infty}=\{u \in H^1_0(D): \text{$u$ is radially symmetric}\}$.
For any $p,q>1$, $\alpha\geq 0$ and $n\in \N\cup\{\infty\}$, set
$$
\San=\inf\{
\Ra(u,v)
:u,v\in H_n\setminus\{0\}
\}.
$$
From the proof of \cite{WY}*{Proposition 2.5},
we can find 
\begin{equation}
\label{Henon-radial}
\Sar
\geq \Szo\Big(\frac{\alpha+2}{2}\Big)^{1+\frac{2}{p+q}}.
\end{equation}
Next,
let $\varphi$ be any element of $C^\infty_0(D)$.
Since we can consider $\varphi\in C^\infty_0(\R^2)$ by the trivial
extension,
we can define $\varphi_\alpha\in C^\infty_0(D)$
by $\varphi_\alpha((x_1,x_2))=\varphi(\alpha (x_1-(1-1/\alpha)), \alpha
x_2)$ for $(x_1,x_2)\in D$.
We set $D_1=D$
and 
$$
D_n=\{(r,\theta): 0< r<1, -\pi/n<\theta< \pi/n\} \quad
\text{for $n \in \N\setminus\{1\}$.}$$
For each $n\in \N$
and $\alpha>0$ with $\mathrm{supp}\,\varphi_\alpha\subset D_n$,
we will show
\begin{equation}
\label{Henon-m-mode}
\San\leq
\Szo n^{1-\frac{2}{p+q}}
\alpha^\frac{4}{p+q}
\Big(\frac{\alpha}{\alpha-2}\Big)^\frac{2\alpha}{p+q}.
\end{equation}
%
We define $P_n: D\rightarrow D$ by
$P_n(r,\theta)=(r,\theta+2\pi/n)$
for $(r,\theta)\in [0,1)\times \R$.
We set $\widetilde{\varphi}(x)=
\varphi_\alpha(x)+\varphi_\alpha(P_n(x))+
\cdots+\varphi_\alpha(P_n^{n-1}(x))$ for $x \in D$.
Since we have
$$\int_D|\nabla \varphi_\alpha|^2=
\int_D|\nabla \varphi|^2
$$
and
$$
\int_D|x|^\alpha|\varphi_\alpha|^p|\varphi_\alpha|^q 
\geq \alpha^{-2}\left(1-\frac{2}{\alpha}\right)^\alpha 
\int_D|\varphi|^p|\varphi|^q,$$
we obtain
$$\San\leq \Ra(\widetilde{\varphi},\widetilde{\varphi})
\leq
\frac{n\displaystyle\int_D(|\nabla\varphi|^2+|\nabla\varphi|^2)}
{\Big(n\alpha^{-2}\Big(1-\frac{2}{\alpha}\Big)^\alpha 
\displaystyle\int_D|\varphi|^p|\varphi|^q\Big)^\frac{2}{p+q}}.
$$
Since $\varphi\in C^\infty_0(D)$ is arbitrary,
we have shown \eqref{Henon-m-mode}.
From \eqref{Henon-radial} and \eqref{Henon-m-mode},
we can see that 
$n \leq n_\alpha$
is a sufficient condition for $\San< \Sar$.
We will show that if
$n>1$
and $\San< \Sar$
then
$\Sao< \cdots<\San$.
Let $n>1$
and $\San< \Sar$.
We can choose $u,v \in H_n\setminus\{0\}$ such that $\Ra(u,v)=\San$
and $u,v\geq 0$.
We note that $u,v \not\in H_\infty$ and $(u,v)$ is a 
positive solution of \eqref{eq:1.1}.
Let $m \in \{1,\ldots,n-1\}$.
We define $\bar{u},\bar{v}\in H_m$ by
$\bar{u}(r,\theta)=u(r, m\theta/n)$
and $\bar{v}(r,\theta)=v(r, m\theta/n)$
for $(r,\theta)\in [0,1)\times\R$.
Since we can see 
$$
\int_D|x|^\alpha |\bar{u}|^p|\bar{v}|^q
=\int_D|x|^\alpha |u|^p|v|^q,
$$
$$
\int_D|\nabla \bar{u}|^2
=\int_0^{2\pi}\int_0^1\Big(
\left|\frac{\partial u}{\partial r}\right|^2 +
\frac{m^2}{n^2r^2} \left|\frac{\partial u}{\partial \theta}\right|^2
\Big)r\,drd\theta< \int_D|\nabla u|^2
$$
and
$\int_D|\nabla \bar{v}|^2< \int_D|\nabla v|^2$,
we have 
$$\Sam\leq \Ra(\bar{u},\bar{v})< \Ra(u,v)=\San.$$
By a similar argument, 
we conclude that
$\Sao< \cdots<\San$.
Hence we infer that if $\San< \Sar$, then
for each $\ell=1,\ldots,n$,
there exists a nonradial positive solution $(u_\ell,v_\ell) \in
H_\ell\times H_\ell$ 
of \eqref{eq:1.1} satisfying
$\Ra(u_\ell,v_\ell)=S_{\alpha,p,q,\ell}$.
We set the number in \eqref{m-alpha} as $\eta(\alpha,p,q)$.
For a fixed $\alpha\in (2,\infty)$,
we have $\eta(\alpha,p,q)\rightarrow 1+\alpha/2$ as $p+q\rightarrow \infty$,
which yields \thetag{i}.
For a fixed $p,q \in (1,\infty)$,
we have $\eta(\alpha,p,q)\rightarrow \infty$ as $\alpha\rightarrow
\infty$, yielding \thetag{ii}.
Hence, we finish the proof of Theorem~\ref{applicationHenon}.

\bigskip
\end{document}